\documentclass[10pt]{article}

\usepackage[all]{xy}
\usepackage{textcomp}
\usepackage{amsmath,amssymb,eucal,eufrak}
\usepackage{graphicx}
\newtheorem{theorem}{Theorem}[section]

\newenvironment{proof}{\par{\bf Proof\ }}{}

\newtheorem{definition}[theorem]{Definition}

\numberwithin{equation}{section}

 \def\NN{{\mathbb{N}}}

  \def\cL{{\mathcal{L}}}

\begin{document}

\title{On concrete models for local operator spaces}

\author{Xavier Mary \footnote{email: xavier.mary@ensae.fr}\\
\textit{\small Ensae - CREST, 3, avenue Pierre Larousse 92245 Malakoff Cedex, France}}

\date{}
\maketitle
\begin{abstract}
In this short note, we propose a concrete analogue of the space $\cL(H)$ for local operator spaces, the multinormed $C^*$-algebra $\displaystyle\prod_{\alpha} \cL(H_{\alpha})$.
\end{abstract}


In this short note, we propose a concrete analogue of the space $\cL(H)$ for local operator spaces. Recall that an operator space can either be defined concretely as a closed subspace
of the the space $\cL(H)$, with the subordinate operator space structure, or abstractly by a Banach space together with a matrix norm (Ruan's theorem \cite{Ruan2000}). An abstract
local operator space is simply a locally convex space with a collection of (separated) matrix seminorms (Effros, Webster \cite{Effros96}). In other words, it is a projective limit of
(abstract) operator spaces. A concrete model for local operator spaces is given in \cite{Dosiev2008}, where the concrete model is an Arens-Michael (complete locally multiplicative)
algebra of unbounded operators on a Hilbert space. We give here another concrete model, the space $$\displaystyle\prod_{\alpha} \cL(H_{\alpha})=\underset{\leftarrow}{\lim}
\cL(H_{\alpha})$$ where the projective limit is taken in the 'operator space' sense \textit{i.e.} the matrix seminorms are given by
$$\forall \alpha\in \Lambda, \, \forall n\in \NN,\: \rho^{(n)}_{\alpha}\left[\left(T_{i,j}\right)\right]=||\left(\left(T_{i,j}\right)_{\alpha}\right)||_{M_n(H_{\alpha})}$$
Remark that this space endowed with coordinatewise multiplication and involution is a multinormed $C^*$-algebra. It is exactly the $*$-algebra of noncommutative continuous functions on
a quantized domain studied by Dosiev \cite{Dosiev2008}.

\begin{definition}[concrete local operator space]
Let $\{H_{\alpha},\alpha\in \Lambda\}$ be a collection of Hilbert spaces. A concrete local operator space is a subspace of $\displaystyle\prod_{\alpha} \cL(H_{\alpha})$ with its
natural local operator space structure.
\end{definition}

\begin{theorem}
Let $E=\underset{\leftarrow_{\alpha}}{\lim} E_{\alpha}$ be an abstract local operator space. Then there exists a realization of $E$ as a concrete local operator space.
\end{theorem}

\begin{proof}
The simple idea is to define subspaces with the induced topology as projective limits. Let $E=\underset{\leftarrow_{\alpha}}{\lim} E_{\alpha}$, and define $\{H_{\alpha},\alpha\in
\Lambda\}$ a collection of Hilbert spaces, such that for all $\alpha \in \Lambda$, $E_{\alpha}$ embeds in $\cL(H_{\alpha})$ completely isometrically (such a collection exists by Ruan's
theorem). Then with have the operator space equality $E_{\alpha}=\underset{\leftarrow}{\lim} \cL(H_{\alpha})$, with projection $j_{\alpha} : E_{\alpha} \rightarrow \cL(H_{\alpha})$ the
canonical injection. By transitivity of projective limits (the transitivity of the projective limits in the 'operator space' sense follows naturally form the transitivity of the
'classical' projective limits), we get
$$E=\underset{\leftarrow}{\lim} E_{\alpha}=\underset{\leftarrow_{\alpha}}{\lim}\underset{\leftarrow}{\lim} \cL(H_{\alpha})= \underset{\leftarrow_{\alpha}}{\lim} \cL(H_{\alpha})$$
and finally $E$ is a subspace of $\prod_{\alpha} \cL(H_{\alpha})=\underset{\leftarrow_{\alpha}}{\lim} \cL(H_{\alpha})$ with the induced local operator space structure..
\end{proof}


\end{document}